\documentclass[10pt]{article}

\usepackage{amsmath}
\usepackage{amsthm}

\newtheorem{theorem}{Theorem}
\newtheorem{definition}[theorem]{Definition}
\newtheorem{assumption}[theorem]{Assumption}
\newtheorem{example}[theorem]{Example}

\newtheorem{remark}[theorem]{Remark}

\newcommand{\MyHeader}[1]{ \begin{flushleft}\textbf{#1}\end{flushleft} }

\title{The equilibrium equations for the lost sales (r,q) policy with inventory level dependent demand and multiple outstanding orders}

\author{
Bryan Johnston \\ 
\emph{ 837 Princeton Rd., Berkley, MI 48072 } \\
\emph{ btj@umich.edu }
}

\date{July 13, 2009}

\begin{document}

\maketitle

\begin{abstract}
We develop the equilibrium equations for the a model generalizing the continuous review (r, q) lost-sales model with constant lead time, multiple outstanding orders and unit Poisson demand. Demand rate is allowed to depend on inventory level. The model is solved for the case r=2,q=2 as well as the general case q=1. An application is given for a system involving lateral transshipments. \\
\\
\emph{Keywords: Inventory, transshipments, Markov processes, multiple orders, lost sales}
\end{abstract}

\section{Introduction}

We are motivated by the situation where the rate at which inventory level at a store changes depends on the inventory level. Such situations arise for example when the amount of product on display affects the rate at which it sells, the existence of lateral transshipments in a multiple store system, and the situation where a fraction of demand is backordered and the rest become lost sales. The reader may see \cite{Urban1} for a survey of models with inventory level dependent demand. We are also motivated by the need \cite{BT1} and \cite{Hill1} to estimate average inventory levels for the usual (r,q) lost sales model with multiple orders outstanding. To these ends we develop a model which generalizes the usual continuous review (r,q) model with constant lead time and multiple orders outstanding. We then derive a system of integral equations for an embedded process, which when solved provides the equilibrium distribution of inventory levels for the system. Even when customer demand at a store does not depend on inventory level, lateral transshipments and partial backlogging can make the effective rate at which inventory level changes depend on the level itself. Thus, to obtain costs for such a system, it is useful to be able to compute the average fraction of time spent at a given inventory level. 

\section{Model}

In this model, we use the phrase \emph{inventory level} to refer to a state variable $l$ which takes integer values. We interpret $l \ge 0$ as a retailer's on-hand inventory (not including pending orders), and $l < 0$ as units which are committed to backorders. We use \emph{effective rate}  to refer to the rate parameter $\lambda_l$ of the exponential random variable $X_l$ with density function $\lambda_l e^{-\lambda_l t}$. We intepret $X_l$ as the amount of time to transition from state $l$ to state $l-1$ in the absence of outstanding orders.

\begin{example}
For a standard lost sales model, we would have $\lambda_l = \lambda$ for $l>0$ and $\lambda_0 = 0$.
\end{example}
\begin{example}
Customers at a retailer demand units one at a time with Poisson rate $\lambda$. When $l \ge C > 0$, lateral transhipments contribute a Poisson demand with rate $\beta$. When $K<l \le 0$, half of the customer demand is backordered and the other half are lost sales. When $l \le K$, no more backorders are allowed by the retailer. In this situation we have
  \begin{equation}
    \lambda_l =
      \left\{
      \begin{array}{rl}
          \lambda + \beta      & \text{ if } l \ge C \\
          \lambda              & \text{ if } 0< l < C \\
          \frac{\lambda}{2}    & \text{ if } K< l \le 0 \\
          0                    & \text{ if } l \le K
      \end{array}
      \right.
  \end{equation}

\end{example}

Replenishment orders for quantity $q$ with constant lead time $\tau$ are placed when the inventory level transitions from $r-kq +1$ to $r-kq$, for $k \ge 0$. This is equivalent to triggering an order when the net inventory position \cite{HW1} transitions from $r+1$ to  $r$. When an order arrives, the inventory level transitions from $l$ to $l+q$.
\begin{assumption}
We require that $\lambda_{l_L} = 0$ for some inventory level $l_L$ and all $l \le l_L$.
\end{assumption}

With this assumption, the maximum number of orders outstanding at any one time is the largest integer $N_0$ such that $N_0 \le \frac{r - l_L}{q} + 1$.

\section{Procedure}

We wish to find the average long run fraction of time $a(l)$ spent in inventory level $l$, and proceed as follows. First, in the spirit of \cite{JT1}, we find an embedded Markov process $M$ with states $s$. State changes in the embedded process are triggered when orders arrive or when orders are placed. Next we calculate the expected amount of time spent at inventory level $l$ while visiting an $M$ state $s$. Then, the state transition probabilities and equilibrium equations of the embedded process $M$ are derived. The equilibrium distribution of $M$ and the law of total expectation are then used to find the fraction of time spent in an inventory level $l$.

\section{Embedded process $M$}

The embedded Markov process $M$ has states $s = (l, t_1, t_2, \ldots , t_{N_0})$ where $l$ is inventory level, and the $t_i$, with $0 \le t_1 \le t_2 \ldots \le t_{N_0} \le \tau$, are the times until oustanding orders are due to arrive. We will use $\vec{0}^k$ to indicate a vector of $k$ zeros. An $M$ state in which there are k outstanding orders can be denoted by
\begin{eqnarray}
s =&  ( l, 0,0, \ldots , 0, t_{N_0 - k +1} , \ldots , t_{N_0}) &   \nonumber \\
  =&  ( l, \vec{0}^{N_0 - k} , t_{N_0 - k +1} , \ldots , t_{N_0}) \nonumber
\end{eqnarray}
where $ l = r-(k-1)q - d $ with $0 \le d < q$ and $0 < t_i \le \tau$.

\subsection{$M$ state transitions}

State transitions of $M$ occur when an order is placed or an order arrives.

\MyHeader{Case 1 - No orders outstanding}

There are no orders outstanding when $s = (l, \vec{0}^{N_0})$ with $r < l \le r+q$. The state $s$ can only transition to state $(r, \vec{0}^{N_0 -1}, \tau)$.

\MyHeader{Case 2 - $k$ orders outstanding, $1 \le k < N_0$}

Let $ s= ( l, \vec{0}^{N_0 - k} , t_{N_0 - k +1} , \ldots , t_{N_0}) $ be the current $M$ state, with \\
$r-kq < l \le r-(k-1)q$. The following can occur :

 \begin{enumerate}
   \item The next reorder point $r-kq$ is hit in a time $t <t_{N_0 - k +1}$ from now. In this case, another order is placed. $M$ transitions to state 
   $$ s'= ( r-kq, \vec{0}^{N_0 - k - 1} , t_{N_0 - k +1} - t , \ldots , t_{N_0}-t, \tau). $$
   \item The order due to arrive in time $t_{N_0 - k +1}$ arrives before level $r-kq$ is hit. Let $d$ with $d< l - r + kq$ be the decrease in inventory level during time $t_{N_0 - k +1}$. 
   $M$ transitions to 
   $$ s'= ( l-d+q, \vec{0}^{N_0 - k + 1} , t_{N_0 - k +2} - t_{N_0 - k +1} , \ldots , t_{N_0}-t_{N_0 - k +1}) $$
   resulting in one fewer outstainding orders.
 \end{enumerate}

\MyHeader{Case 3 - $N_0$ orders outstanding}

The current state is $s = (l, t_1, t_2, \ldots , t_{N_0})$ with \\
$r-N_0q < l_L \le l \le r-(N_0 - 1)q$. The next order is due to arrive in time $t_1$. Since $\lambda_{l_L} =0$ the inventory level can only decrease by $d \le l-l_L$ during time $t_1$. $M$ will transition to a state
$$s = (l-d+q, 0, t_2 - t_1, \ldots , t_{N_0} - t_1).$$

\begin{remark}
A) It is clear that the only $M$ states with $N_0$ outstanding orders are $s = (r-(N_0 - 1)q, t_1, t_2, \ldots , t_{N_0-1},\tau)$.\\
B) When the system hits a reorder point from above, the next $M$ state has $t_{N_0} = \tau $. \\
C) When an order arrives, the next $M$ state has $t_{N_0} \le \tau $, and $t_{N_0}$ is almost surely strictly less than $\tau$ ( $t_{N_0} = \tau$ requires new order to be placed instantly before a pending order arrives).
\end{remark}

\subsection{Expected time at inventory levels while in an $M$ state}

Let $c(s,l)$ denote the expected time spent at inventory level $l$ while visiting $M$ state $s$. 

\MyHeader{Case 1 - No orders outstanding}

The $M$ state vistited is $s = (l', \vec{0}^{N_0})$ where $r < l' \le r+q$. Only inventory levels $l$ with $r <l \le l'$ can be hit while in $M$ state $s$. There are no outstanding orders, so if inventory level $l$ is hit, the system remains at level $l$ for an expected time $\frac{1}{\lambda_l}$. So we have for $s = (l', \vec{0}^{N_0})$ with $r<l' \le r+q$

\begin{equation}
    c(s,l) =
      \left\{
      \begin{array}{rl}
          \frac{1}{\lambda_l}  & \text{ if } r <l \le l'\\
          0                    & \text{ otherwise }
      \end{array}
      \right.
\end{equation}

\MyHeader{Case 2 - $k$ orders outstanding, $1 \le k < N_0$}

Let $ s= ( l', \vec{0}^{N_0 - k} , t_{N_0 - k +1} , \ldots , t_{N_0}) $ be the current $M$ state, with \\
$r-kq < l' \le r-(k-1)q$. From this state, the system can hit inventory levels $l$ with $r-kq < l \le l'$. The time until the next order arrives is $t_{N_0 - k +1}$ so the expected time spent in level $l$ while in $M$ state $s$ 
is given by \\
$h(l',l,t_{N_0 - k +1})$, where $h()$ is given in the appendix. \\
So we have for $ s= ( l', \vec{0}^{N_0 - k} , t_{N_0 - k +1} , \ldots , t_{N_0}) $ with $r-kq < l' \le r-(k-1)q$

\begin{equation}
    c(s,l) =
      \left\{
      \begin{array}{rl}
         h(l',l,t_{N_0 - k +1})  & \text{ if } r-kq < l \le l'\\
          0                    & \text{ otherwise }
      \end{array}
      \right.
\end{equation}

\MyHeader{Case 3 - $N_0$ orders outstanding}

The curent state is $s = (l', t_1, t_2, \ldots , t_{N_0})$ with \\
$r-N_0q < l_L \le l' = r-(N_0 - 1)q$. The next order is due to arrive in time $t_1$. Similar to Case 2 we have

\begin{equation}
    c(s,l) =
      \left\{
      \begin{array}{rl}
         h(r-(N_0 - 1)q,l,t_1)  & \text{ if } r-N_0q < l_L \le l \le l' = r-(N_0 - 1)q\\
          0                    & \text{ otherwise }
      \end{array}
      \right.
\end{equation}

\subsection{$M$ state transition probabilities}
Let $p(s_1,s_2)$ denote the probability mass or probability density of transitioning from $M$ state $s_1$ to $M$ state $s_2$. 

\MyHeader{Case 1 - No orders outstanding}

The current $M$ state vistited is $s_1 = (l_1, \vec{0}^{N_0})$ where $r < l_1 \le r+q$. Clearly

\begin{equation}
    p(s_1,s_2) =
      \left\{
      \begin{array}{rl}
         1  & \text{ if } s_2 = (r,\vec{0}^{N_0 -1}, \tau) \\
         0  & \text{ otherwise }
      \end{array}
      \right.
\end{equation}

\MyHeader{Case 2 - $k$ orders outstanding, $1 \le k < N_0$}

Let $ s_1 = ( l_1, \vec{0}^{N_0 - k} , t_{N_0 - k +1} , \ldots , t_{N_0}) $ be the current $M$ state, with \\
$r-kq < l_1 \le r-(k-1)q$. The next order is due to arrive in time $t_{N_0 - k +1}$. We saw in section 4.1, case 2, that M will transition to one of 

\begin{enumerate}
 \item $ s_2= ( r-kq, \vec{0}^{N_0 - k - 1} , t_{N_0 - k +1} - t , \ldots , t_{N_0}-t, \tau)$  , if the inventory level decreases by $l_1 - r + kq$ in time $t \le t_{N_0 - k +1}$. This occurs with probability density $f(l_1, r-kq +1, t)$ (see appendix).
 \item $ s_2= ( l_1-d+q, \vec{0}^{N_0 - k + 1} , t_{N_0 - k +2} - t_{N_0 - k +1} , \ldots , t_{N_0}-t_{N_0 - k +1})$ if the inventory level decreases in time $t_{N_0 - k +1}$ by $d$ with $d < l_1 - r + kq$. 
 This occurs with probability mass $g(l_1,d, t_{N_0 - k +1})$ (see appendix).
\end{enumerate}

So $p(s_1,s_2)$ is a mixed discrete continuous distribution, and is given by\\
\\
$p(s_1,s_2) = $ 
\begin{equation}
    \left\{
    \begin{array}{rl}
       f(l_1, r-kq +1, t)  
            & \text{ if } s_2= ( r-kq, \vec{0}^{N_0 - k - 1} , t_{N_0 - k +1} - t , \ldots , t_{N_0}-t, \tau)\\
            &  \text{ with } t \le t_{N_0 - k +1}\\
       g(l_1,d, t_{N_0 - k +1})
            & \text{ if } s_2= ( l_1-d+q, \vec{0}^{N_0 - k + 1} , t_{N_0 - k +2} - t_{N_0 - k +1} , \ldots , t_{N_0}-t_{N_0 - k +1}) \\
            &  \text{ with } d < l_1 - r + kq \\
       0  & \text{ otherwise }
    \end{array}
    \right.
\end{equation}

\MyHeader{Case 3 - $N_0$ orders outstanding}

The current state is $s_1 = (l_1, t_1, t_2, \ldots , t_{N_0})$ with \\
$r-N_0q < l_L \le l_1 = r-(N_0 - 1)q$. The next order is due to arrive in time $t_1$. Since $\lambda_{l_L} = 0$, this $M$ state will always transition to one with one fewer outstanding orders. So we have

$p(s_1,s_2) = $ 
\begin{equation}
    \left\{
    \begin{array}{rl}
       g(l',d, t_{N_0 - k +1})
            & \text{ if } s_2= (l'- d + q, 0, t_2 - t_1, \ldots , t_{N_0} - t_1) \\
            &  \text{ with } d \le l-l_L \text{ and } l' = r-(N_0 - 1)q\\
       0  & \text{ otherwise }
    \end{array}
    \right.
\end{equation}

\section{Equilibrium distribution for embedded process $M$}

The long-run distribution of the states of the embedded process $M$ will have a mixed discrete continuous distribution. It will be convenient to use different symbols to represent contributions to the distribution.

\begin{definition} For the embedded process $M$ we define \\
A) \emph{(Order arrives)} $\bar{y}(l, \vec{0}^{N_0})$ for $r<l \le r+q$ is the long-run probability mass that $M$ transitions to the state $(l, \vec{0}^{N_0})$ with no outstanding orders. \\
B) \emph{(New order placed)} $\hat{y}(r, \vec{0}^{N_0 - 1}, \tau)$ is the long-run probability mass that $M$
transitions to state $(r, \vec{0}^{N_0 - 1}, \tau)$. \\
C) \emph{(New order placed)} $\hat{y}(r-(k-1)q, \vec{0}^{N_0 - k} , t_{N_0 - k + 1}, \ldots , t_{N_0 - 1}, \tau)$, \\
for $2 \le k < N_0$, is the long-run probability density that $M$ \\
transitions to state $(r-(k-1)q, \vec{0}^{N_0 - k} , t_{N_0 - k + 1}, \ldots , t_{N_0 - 1}, \tau)$. \\
D) \emph{(Order arrives)} $y(l, \vec{0}^{N_0 - k} , t_{N_0 - k + 1}, \ldots , t_{N_0 - 1}, t_{N_0})$, \\
for $1 \le k \le N_0 - 1$, $r-kq < l \le r - (k-1)q$ and $t_{N_0} < \tau$ ,\\
is the long-run probability density that $M$ \\
transitions to state $(l, \vec{0}^{N_0 - k} , t_{N_0 - k + 1}, \ldots , t_{N_0 - 1}, t_{N_0})$.
\end{definition}

We see that for fixed $t_1, \ldots, t_{N_0 -1}$, the distribution is continuous for $t_{N_0} < \tau$ 
and discrete at $t_{N_0} = \tau$.

We will proceed on a case-by-case basis to write the equilibrium euqations for the embedded process $M$.

\MyHeader{Case 1 - $\bar{y}(l, \vec{0}^{N_0})$ }

$M$ state $s = (l, \vec{0}^{N_0})$ , $r < l \le r+q$, can only be hit from states $s' = (r, \vec{0}^{N_0 - 1}, t_{N_0})$ with $t_{N_0} \le \tau $. We use the transition probabilites from section 4.3, integrate the continuous contributions and add the discrete contributions to obtain

\begin{eqnarray*}
  \bar{y}(l, \vec{0}^{N_0}) & = & \sum_{l' = l-q}^{r}{
                                         \int_{x=0}^{\tau}{
                                          p((l', \vec{0}^{N_0 - 1}, t_{N_0}), (l, \vec{0}^{N_0}) )
                                          y(l', \vec{0}^{N_0 - 1}, x) dx 
                                         }
                                        } \\
                                      && + p((r, \vec{0}^{N_0 - 1}, \tau), (l, \vec{0}^{N_0}) )
                                      \hat{y}(r, \vec{0}^{N_0 - 1}, \tau) \\
                           & = & \sum_{l' = l-q}^{r}{
                                         \int_{x=0}^{\tau}{
                                          g(l',l' + q -l, x) y(l', \vec{0}^{N_0 - 1}, x) dx 
                                         }
                                        } \\
                                      && + g(r,r + q - l, \tau) \hat{y}(r, \vec{0}^{N_0 - 1}, \tau)
\end{eqnarray*}

\MyHeader{Case 2 - $\hat{y}(r, \vec{0}^{N_0 - 1}, \tau)$ }

$M$ state $s = (r, \vec{0}^{N_0 - 1}, \tau)$ can only be hit from states $s' = (l, \vec{0}^{N_0})$ , $r < l \le r+q$. Each of these states transitions to state $(r, \vec{0}^{N_0 - 1}, \tau)$ with probability 1, so

\begin{eqnarray*}
  \hat{y}(r, \vec{0}^{N_0 - 1}, \tau) &=& \sum_{l' = r+1}^{r+q}{ \bar{y}(l', \vec{0}^{N_0}) }
\end{eqnarray*}

\MyHeader{Case 3 - $\hat{y}(r-(k-1)q, \vec{0}^{N_0 - k} , t_{N_0 - k + 1}, \ldots , t_{N_0 - 1}, \tau)$, for $2 \le k < N_0$ }

$M$ state $s=(r-(k-1)q, \vec{0}^{N_0 - k} , t_{N_0 - k + 1}, \ldots , t_{N_0 - 1}, \tau)$, for $2 \le k < N_0$, can only be hit from states $s'=(l, \vec{0}^{N_0 - k +1} , t + t_{N_0 - k + 1}, \ldots ,t +  t_{N_0 - 1})$ with $t_{N_0 - 1} \le \tau$ and $t \le \tau - t_{N_0 - 1}$. This occurs when a new order is placed before an outstanding order arrives. Summing and integrating contributions gives

\begin{flushleft}$\hat{y}(r-(k-1)q, \vec{0}^{N_0 - k} , t_{N_0 - k + 1}, \ldots , t_{N_0 - 1}, \tau) $\end{flushleft}
\begin{eqnarray*}
&& = \\
&&\sum_{l' = r-(k-1)q+1}^{r-(k-2)q}{} \Bigg[ \int_{x=0}^{\tau}{} \bigg\{ \\
  && \qquad \qquad p((l', \vec{0}^{N_0 - k +1} , x + t_{N_0 - k + 1}, \ldots ,x +  t_{N_0 - 1}), s ) \\
  && \qquad \qquad * \quad y(l', \vec{0}^{N_0 - k +1} , x + t_{N_0 - k + 1}, \ldots ,x +  t_{N_0 - 1}) \\
  && \qquad \bigg\} dx  \Bigg] \\
  &&  + \\
  &&  \qquad  p((r-(k-2)q, \vec{0}^{N_0 - k +1} ,\\
  && \qquad \qquad  (\tau - t_{N_0 - 1}) + t_{N_0 - k + 1}, \ldots ,(\tau - t_{N_0 - 1}) +  t_{N_0 - 2}, \tau), s ) \\
  && \qquad  * \quad \hat{y}(r-(k-2)q, \vec{0}^{N_0 - k +1} , \\
  && \qquad  \qquad \quad (\tau - t_{N_0 - 1}) + t_{N_0 - k + 1}, \ldots ,(\tau - t_{N_0 - 1}) +  t_{N_0 - 2}, \tau) \\
&& = \\
  && \sum_{l' = r-(k-1)q+1}^{r-(k-2)q}{} \Bigg[ \int_{x=0}^{\tau}{} \bigg\{ \\
  && \qquad \qquad  f(l', r-(k-1)q+1, x) * y(l', \vec{0}^{N_0 - k +1} , x + t_{N_0 - k + 1}, \ldots ,x +  t_{N_0 - 1}) \\
  && \qquad \bigg\} dx  \Bigg] \\
  && + \\
  && \qquad f(r-(k-2)q, r-(k-1)q+1, \tau - t_{N_0 - 1}) \\
  && \qquad * \quad \hat{y}(r-(k-2)q, \vec{0}^{N_0 - k +1} , \\
  && \qquad \qquad  (\tau - t_{N_0 - 1}) + t_{N_0 - k + 1}, \ldots ,(\tau - t_{N_0 - 1}) +  t_{N_0 - 2}, \tau)
\end{eqnarray*}

\MyHeader{Case 4 - $y(l, \vec{0}^{N_0 - k} , t_{N_0 - k + 1}, \ldots , t_{N_0 - 1}, t_{N_0})$, for $1 \le k \le N_0 - 2$ }

$M$ state $s = (l, \vec{0}^{N_0 - k} , t_{N_0 - k + 1}, \ldots , t_{N_0 - 1}, t_{N_0})$, for $1 \le k \le N_0 - 2$, $t_{N_0} < \tau$ and $r-kq < l \le r-(k-1)q$, can only be hit from $M$ states $s' = (l', \vec{0}^{N_0 - k - 1} , t, t + t_{N_0 - k + 1}, \ldots , t + t_{N_0 - 1}, t + t_{N_0})$, with $l-q \le l' \le r-kq$ and $0 < t \le \tau-t_{N_0}$. This occurs when the order pending arrival in time $t$ arrives before the inventory level decreases to the next reorder point at $r-(k+1)q$. As above, summing and integrating contributions gives

\begin{flushleft}$y(l, \vec{0}^{N_0 - k} , t_{N_0 - k + 1}, \ldots , t_{N_0 - 1}, t_{N_0})$\end{flushleft}
\begin{eqnarray*}
&& = \\
&&\sum_{l' = l-q}^{r-kq}{} \Bigg[ \int_{x=0}^{\tau - t_{N_0}}{} \bigg\{ \\
  && \qquad \qquad p((l', \vec{0}^{N_0 - k - 1} , t, t + t_{N_0 - k + 1}, \ldots , t + t_{N_0 - 1}, t + t_{N_0}), s ) \\
  && \qquad \qquad * \quad y(l', \vec{0}^{N_0 - k - 1} , t, t + t_{N_0 - k + 1}, \ldots , t + t_{N_0 - 1}, t + t_{N_0}) \\
  && \qquad \bigg\} dx  \Bigg] \\
  &&  + \\
  &&  \qquad  p((r-kq, \vec{0}^{N_0 - k - 1} , (\tau-t_{N_0}), (\tau-t_{N_0}) + t_{N_0 - k + 1}, \ldots , (\tau-t_{N_0}) + t_{N_0 - 1}, \tau), s ) \\
  && \qquad  * \quad \hat{y}(r-kq, \vec{0}^{N_0 - k - 1} , (\tau-t_{N_0}), (\tau-t_{N_0}) + t_{N_0 - k + 1}, \ldots , (\tau-t_{N_0}) + t_{N_0 - 1}, \tau) \\
&& = \\
  && \sum_{l' = l-q}^{r-kq}{} \Bigg[ \int_{x=0}^{\tau - t_{N_0}}{} \bigg\{ \\
  && \qquad \qquad  g(l', l'+q-l, t) * y(l', \vec{0}^{N_0 - k - 1} , t, t + t_{N_0 - k + 1}, \ldots , t + t_{N_0 - 1}, t + t_{N_0}) \\
  && \qquad \bigg\} dx  \Bigg] \\
  && + \\
  && \qquad g(r-kq, r-(k-1)q -l,\tau-t_{N_0} ) \\
  && \qquad * \quad \hat{y}(r-kq, \vec{0}^{N_0 - k - 1} , (\tau-t_{N_0}), (\tau-t_{N_0}) + t_{N_0 - k + 1}, \ldots , (\tau-t_{N_0}) + t_{N_0 - 1}, \tau)
\end{eqnarray*}

\MyHeader{Case 5 - $y(l, \vec{0}^{N_0 - k} , t_{N_0 - k + 1}, \ldots , t_{N_0 - 1}, t_{N_0})$, for $k = N_0 - 1$ }

These are states of the form $s = (l, 0 , t_2, \ldots , t_{N_0 - 1}, t_{N_0})$, with $t_{N_0} < \tau$ and $r-(N_0 - 1)q < l \le r-N_0 - 2)q$. We observed before that there are no $M$ states with $N_0$ outstanding orders with inventory level less that $r-(N_0-1)q$. So the states $s$ can only be hit (when a pending order arrives) from $M$ states $s' = (r-(N_0-1)q , (\tau-t_{N_0}), (\tau-t_{N_0}) + t_2, \ldots , (\tau-t_{N_0}) + t_{N_0 - 1}, (\tau-t_{N_0}) + t_{N_0} = \tau)$. Similar to the previous case we find

\begin{flushleft}$y(l, 0 , t_2, \ldots , t_{N_0 - 1}, t_{N_0})$\end{flushleft}
\begin{eqnarray*}
&& = \\
 &&  \qquad  p((r-(N_0-1)q , \tau-t_{N_0}, (\tau-t_{N_0}) + t_2, \ldots , (\tau-t_{N_0}) + t_{N_0 - 1}, \tau), s ) \\
  && \qquad  * \quad \hat{y}(r-(N_0-1)q , \tau-t_{N_0}, (\tau-t_{N_0}) + t_2, \ldots , (\tau-t_{N_0}) + t_{N_0 - 1},  \tau) \\
&& = \\
  && \qquad g(r-(N_0-1)q, r-(N_0-2)q - l,\tau-t_{N_0} ) \\
  && \qquad * \quad \hat{y}(r-(N_0-1)q , \tau-t_{N_0}, (\tau-t_{N_0}) + t_2, \ldots , (\tau-t_{N_0}) + t_{N_0 - 1}, \tau)
\end{eqnarray*}

\MyHeader{Normalization equation}

The above equations have solutions $y()$, $\bar{y}()$ and $\hat{y}()$ which are only unique up to a constant multiple. We are concerned with the average fraction of time $a(l)$ spent at an inventory level $l$. Solving the equilibrium equations for $M$ up to a constant multiple will suffice for that purpose. For completeness, however, we present the normalization equation for the embedded process. By integrating contributions from $y()$, $\bar{y}()$ and $\hat{y}()$ over the state space of $M$ we have

\begin{eqnarray*}
1 = && \\
    && \sum_{l' = r+1}^{r+q}{ \bar{y}(l,\vec{0}^{N_0}) } \\
    &+& \hat{y}(r, \vec{0}^{N_0 - 1}, \tau) \\
    &+& \sum_{k=1}^{N_0-1}{} \sum_{l = r-kq+1}^{r-(k-1)q}{} \Bigg\{ 
        \int_{t_{N_0} =0}^{\tau}{}  \int_{t_{N_0 -1} =0}^{t_{N_0}}{} \ldots  
              \int_{t_{N_0 - k +1} =0}^{t_{N_0 -k +2}}{} \\
      &&  \qquad \qquad \qquad \qquad y(l, \vec{0}^{N_0 - k} , t_{N_0 - k + 1}, \ldots , t_{N_0 - 1}, t_{N_0}) 
          dt_{N_0 - k +1} \ldots dt_{N_0} \Bigg\} \\
    &+& \sum_{k=2}^{N_0}{} \Bigg\{ 
       \int_{t_{N_0-1} =0}^{\tau}{}  \int_{t_{N_0 -2} =0}^{t_{N_0-1}}{} \ldots 
            \int_{t_{N_0 - k +1} =0}^{t_{N_0 -k +2}}{} \\
      && \qquad \qquad \hat{y}(r-(k-1)q, \vec{0}^{N_0 - k} , t_{N_0 - k + 1}, \ldots , t_{N_0 - 1}, \tau)
        dt_{N_0 - k +1} \ldots dt_{N_0 -1} \Bigg\}
\end{eqnarray*}

\subsection{The equilibrium equations for $M$}

We have just shown that the equilibrium equations for the embedded process are

\begin{enumerate}
  \item{$r < l \le r+q$}
    \begin{eqnarray}
      \bar{y}(l, \vec{0}^{N_0}) & = & \sum_{l' = l-q}^{r}{
                                             \int_{x=0}^{\tau}{
                                              g(l',l' + q -l, x) y(l', \vec{0}^{N_0 - 1}, x) dx 
                                             }
                                            } \nonumber \\
                                          && + g(r,r + q - l, \tau) \hat{y}(r, \vec{0}^{N_0 - 1}, \tau)
    \end{eqnarray}

  \item
    \begin{eqnarray}
      \hat{y}(r, \vec{0}^{N_0 - 1}, \tau) &=& \sum_{l' = r+1}^{r+q}{ \bar{y}(l', \vec{0}^{N_0}) } \nonumber \\
      && 
    \end{eqnarray}

  \item{$2 \le k < N_0$} \\
    \\
    $\hat{y}(r-(k-1)q, \vec{0}^{N_0 - k} , t_{N_0 - k + 1}, \ldots , t_{N_0 - 1}, \tau) $
    \begin{eqnarray}
    &=& \nonumber  \\
      && \sum_{l' = r-(k-1)q+1}^{r-(k-2)q}{} \Bigg[ \int_{x=0}^{\tau}{} \bigg\{ \nonumber  \\
      && \qquad \qquad  f(l', r-(k-1)q+1, x) 
          * y(l', \vec{0}^{N_0 - k +1} , x + t_{N_0 - k + 1}, \ldots ,x +  t_{N_0 - 1}) \nonumber  \\
      && \qquad \bigg\} dx  \Bigg] \nonumber  \\
      && + \nonumber  \\
      && \qquad f(r-(k-2)q, r-(k-1)q+1, \tau - t_{N_0 - 1}) \nonumber  \\
      && \qquad * \quad \hat{y}(r-(k-2)q, \vec{0}^{N_0 - k +1} , \nonumber  \\
      && \qquad \qquad  (\tau - t_{N_0 - 1}) + t_{N_0 - k + 1}, \ldots ,(\tau - t_{N_0 - 1}) +  t_{N_0 - 2}, \tau) \nonumber \\
      &&
    \end{eqnarray}

  \item{$1 \le k \le N_0 - 2$, $t_{N_0} < \tau$ and $r-kq < l \le r-(k-1)q$} \\
    \\
    $y(l, \vec{0}^{N_0 - k} , t_{N_0 - k + 1}, \ldots , t_{N_0 - 1}, t_{N_0})$
    \begin{eqnarray}
    &=& \nonumber  \\
      && \sum_{l' = l-q}^{r-kq}{} \Bigg[ \int_{x=0}^{\tau - t_{N_0}}{} \bigg\{ \nonumber  \\
      && \qquad \qquad  g(l', l'+q-l, t) 
          * y(l', \vec{0}^{N_0 - k - 1} , t, t + t_{N_0 - k + 1}, \ldots , t + t_{N_0 - 1}, t + t_{N_0}) \nonumber  \\
      && \qquad \bigg\} dx  \Bigg] \nonumber  \\
      && + \nonumber  \\
      && \qquad g(r-kq, r-(k-1)q -l,\tau-t_{N_0} ) \nonumber  \\
      && \qquad * \quad \hat{y}(r-kq, \vec{0}^{N_0 - k - 1} , (\tau-t_{N_0}), (\tau-t_{N_0}) + t_{N_0 - k + 1}, \ldots , (\tau-t_{N_0}) + t_{N_0 - 1}, \tau) \nonumber \\
      && 
    \end{eqnarray}

  \item{$t_{N_0} < \tau$ and $r-(N_0 - 1)q < l \le r-N_0 - 2)q$} \\
    \\
    $y(l, 0 , t_2, \ldots , t_{N_0 - 1}, t_{N_0})$
    \begin{eqnarray}
    &=& \nonumber  \\
      && \qquad g(r-(N_0-1)q, r-(N_0-2)q - l,\tau-t_{N_0} )\nonumber   \\
      && \qquad * \quad 
          \hat{y}(r-(N_0-1)q , \tau-t_{N_0}, (\tau-t_{N_0}) + t_2, \ldots , (\tau-t_{N_0}) + t_{N_0 - 1}, \tau) \nonumber \\
      && 
    \end{eqnarray}

  \item{Normalization}
    \begin{eqnarray}
    1 = && \nonumber  \\
        && \sum_{l' = r+1}^{r+q}{ \bar{y}(l,\vec{0}^{N_0}) } \nonumber  \\
        &+& \hat{y}(r, \vec{0}^{N_0 - 1}, \tau)\nonumber   \\
        &+& \sum_{k=1}^{N_0-1}{} \sum_{l = r-kq+1}^{r-(k-1)q}{} \Bigg\{ 
            \int_{t_{N_0} =0}^{\tau}{}  \int_{t_{N_0 -1} =0}^{t_{N_0}}{} \ldots  
                  \int_{t_{N_0 - k +1} =0}^{t_{N_0 -k +2}}{} \nonumber  \\
          &&  \qquad \qquad \qquad \qquad y(l, \vec{0}^{N_0 - k} , t_{N_0 - k + 1}, \ldots , t_{N_0 - 1}, t_{N_0}) 
              dt_{N_0 - k +1} \ldots dt_{N_0} \Bigg\} \nonumber  \\
        &+& \sum_{k=2}^{N_0}{} \Bigg\{ 
           \int_{t_{N_0-1} =0}^{\tau}{}  \int_{t_{N_0 -2} =0}^{t_{N_0-1}}{} \ldots 
                \int_{t_{N_0 - k +1} =0}^{t_{N_0 -k +2}}{} \nonumber  \\
          && \qquad \qquad \hat{y}(r-(k-1)q, \vec{0}^{N_0 - k} , t_{N_0 - k + 1}, \ldots , t_{N_0 - 1}, \tau)
            dt_{N_0 - k +1} \ldots dt_{N_0 -1} \Bigg\} \nonumber \\
         &&
    \end{eqnarray}

\end{enumerate}
\begin{flushright} $\qed$ \end{flushright}

\section{Average fraction of time spent at an inventory level}

Let $w(l)$ denote the expected time per $M$ state spent at inventory level $l$. We can compute $w(l)$ using the distribution for the embedded processes $M$ and the law of total expectation. After computing $w(l)$, we can compute the average fraction of time $a(l)$ spent at inventory level $l$ as 

\begin{equation}
  a(l) = \frac{w(l)}{ \sum_{i=l_L}^{r+q}{w(i)} }
\end{equation}

We'll proceed as before and break the situation into cases.

\MyHeader{Case 1 - $r<l \le r+q$}
Inventory level $l$ can only be hit in $M$ states $s = (l', \vec{0}^{N_0})$ with $l \le l' \le r+q$. So applying the law of total expectation gives

\begin{equation}
  w(l) = \sum_{l'=l}^{r+q}{c((l', \vec{0}^{N_0}),l) \bar{y}(l', \vec{0}^{N_0}) }
\end{equation}

\MyHeader{Case 2 - $r-kq < l \le r-(k-1)q$ for $1 \le k \le N_0 - 1$ }
Level $l$ can only be hit in $M$ states 

\begin{enumerate}
  \item
    $s = (l', \vec{0}^{N_0 - k}, t_{N_0 - k + 1}, \ldots , t_{N_0 - 1}, t_{N_0})$ with $l \le \l' \le r-(k-1)q$ and $t_{N_0}<\tau$, and
  \item
    $s = (r-(k-1)q, \vec{0}^{N_0 - k}, t_{N_0 - k + 1}, \ldots , t_{N_0 - 1}, \tau)$
\end{enumerate}

So the total expected time per $M$ state spent at level $l$ is

\begin{eqnarray}
w(l) &=& \sum_{l'=l}^{r-(k-1)q}{} 
         \int_{t_{N_0} =0}^{\tau} \int_{t_{N_0 -1} =0}^{t_{N_0}} 
              \ldots \int_{t_{N_0 -k+1} =0}^{t_{N_0 -k+2}}  \Bigg\{ \nonumber \\
      && \qquad c((l', \vec{0}^{N_0 -k},t_{N_0 - k +1}, \ldots t_{N_0} ),l) \nonumber \\
      &&  \qquad \qquad *  y(l', \vec{0}^{N_0 -k},t_{N_0 - k +1}, \ldots t_{N_0} ) \nonumber \\
      && \qquad \qquad \qquad \qquad \qquad dt_{N_0 - k +1} \ldots dt_{N_0} \Bigg\} \nonumber \\
      &+& \int_{t_{N_0 -1} =0}^{\tau} \int_{t_{N_0 -2} =0}^{t_{N_0 -1}} 
              \ldots \int_{t_{N_0 -k+1} =0}^{t_{N_0 -k+2}}  \Bigg\{ \nonumber \\
      && \qquad c((r-(k-1)q, \vec{0}^{N_0 -k},t_{N_0 - k +1}, \ldots t_{N_0-1}, \tau ),l) \nonumber \\
      &&  \qquad \qquad * \hat{y}(r-(k-1)q, \vec{0}^{N_0 -k},t_{N_0 - k +1}, \ldots t_{N_0-1}, \tau ) \nonumber \\
      && \qquad \qquad \qquad \qquad \qquad dt_{N_0 - k +1} \ldots dt_{N_0 -1} \Bigg\} \nonumber \\
\end{eqnarray}

\MyHeader{Case 3 - $l_L \le l \le r-(N_0-1)q$ for $1 \le k \le N_0 - 1$ }
Level $l$ can only be hit in $M$ states $s = (r-(N_0-1)q, t_1, t_2, \ldots , t_{N_0 - 1}, \tau)$.

\begin{eqnarray}
w(l) &=& \int_{t_{N_0 -1} =0}^{\tau} \int_{t_{N_0 -2} =0}^{t_{N_0 -1}} 
              \ldots \int_{t_1 =0}^{t_2}  \Bigg\{ \nonumber \\
      && \qquad c((r-(N_0-1)q, t_1, t_2, \ldots , t_{N_0 - 1}, \tau),l) \nonumber \\
      &&  \qquad \qquad * \hat{y}(r-(N_0-1)q, t_1, t_2, \ldots , t_{N_0 - 1}, \tau) \nonumber \\
      && \qquad \qquad \qquad \qquad \qquad dt_1 \ldots dt_{N_0 -1} \Bigg\} \nonumber \\
\end{eqnarray}

\section{The case q=1}

In this section we will solve the equilibrium equations for the case $q=1$ and find a closed form expression for $a(l)$. Referring to the appendix we have $f(l,l,x) = \lambda_l e^{-\lambda_l x}$ and $g(l,0,x) = e^{-\lambda_l x}$. To make verifying the solutions more convenient (we will leave out these details), we will rewrite the equilibrium equations of $M$ as 

\begin{enumerate}
  \item{Equation 0A}

  \begin{eqnarray}
    \bar{y}(r+1, \vec{0}^{N_0} )
        &=& \int_{x=0}^{\tau}{ e^{-\lambda_r x} y(r,\vec{0}^{N_0 -1}, x) dx } \nonumber \\
        && \qquad + e^{-\lambda_r \tau} \hat{y}(r,\vec{0}^{N_0 -1}, \tau)
  \end{eqnarray}

  \item{Equation 0B}

  \begin{eqnarray}
    \hat{y}(r,\vec{0}^{N_0 -1}, \tau) &=& \bar{y}(r+1, \vec{0}^{N_0} )
  \end{eqnarray}

  \item{Equation 0C}

  \begin{eqnarray}
    y(r, \vec{0}^{N_0 -1} ,t)
        &=& \int_{x=0}^{\tau - t}{ e^{-\lambda_{r-1} x} y(r-1,\vec{0}^{N_0 -2}, x, x+t) dx } \nonumber \\
        && \qquad + e^{-\lambda_{r-1}(\tau -t)} \hat{y}(r-1,\vec{0}^{N_0 -2}, \tau - t, \tau)
  \end{eqnarray}

  \item{Equation kB, for $2 \le k \le N_0 -2$ }
  
  \begin{flushleft}$\hat{y}(r - (k-1), \vec{0}^{N_0 -k} ,t_{N_0 - k +1}, \ldots, t_{N_0-1}, \tau)$\end{flushleft}
  \begin{eqnarray}
    &=& \int_{x=0}^{\tau - t_{N_0-1}}{} \lambda_{r-(k-2)} e^{-\lambda_{r-(k-2)} x} \nonumber \\
    && \qquad * y(r-(k-2),\vec{0}^{N_0 -k +1}, x + t_{N_0 -k +1}, x + t_{N_0 - 1}) dx  \nonumber \\
    && + \lambda_{r-(k-2)} e^{-\lambda_{r-(k-2)} (\tau - t_{N_0 - 1}) } \nonumber \\
    &&  \qquad * \hat{y}(r-(k-2),\vec{0}^{N_0 -k +1}, 
            (\tau - t_{N_0 - 1}) + t_{N_0 -k +1}, \ldots , (\tau - t_{N_0 - 1}) + t_{N_0 -2}, \tau) \nonumber \\
    &&
  \end{eqnarray}

  \item{Equation kC, for $2 \le k \le N_0 -2$ }
  
  \begin{flushleft}$y(r - (k-1), \vec{0}^{N_0 -k} ,t_{N_0 - k +1}, \ldots ,t_{N_0-1}, t_{N_0})$\end{flushleft}
  \begin{eqnarray}
    &=& \int_{x=0}^{\tau - t_{N_0}}{} e^{-\lambda_{r-k} x} \nonumber \\
    && \qquad * y(r-k,\vec{0}^{N_0 -k -1}, x, x + t_{N_0 -k +1}, \ldots , x + t_{N_0}) dx  \nonumber \\
    && + e^{-\lambda_{r-k} (\tau - t_{N_0}) } \nonumber \\
    &&  \qquad * \hat{y}(r-k,\vec{0}^{N_0 -k -1}, 
            \tau - t_{N_0}, (\tau - t_{N_0}) + t_{N_0 -k +1}, \ldots , (\tau - t_{N_0}) + t_{N_0 -1}, \tau) \nonumber \\
    &&
  \end{eqnarray}

  \item{Equation ($N_0 -1$)B}
  
  \begin{flushleft}$\hat{y}(r - (N_0 -2), 0 ,t_2, \ldots, t_{N_0-1}, \tau)$\end{flushleft}
  \begin{eqnarray}
    &=& \int_{x=0}^{\tau - t_{N_0-1}}{} \lambda_{r-(N_0-3)} e^{-\lambda_{r-(N_0-3)} x} \nonumber \\
    && \qquad * y(r-(N_0-3)),0,0, x + t_2, x + t_{N_0 - 1}) dx  \nonumber \\
    && + \lambda_{r-(N_0-3)} e^{-\lambda_{r-(N_0-3)} (\tau - t_{N_0 - 1}) } \nonumber \\
    &&  \qquad * \hat{y}(r-(N_0-3),0,0, 
            (\tau - t_{N_0 - 1}) + t_2, \ldots , (\tau - t_{N_0 - 1}) + t_{N_0 -2}, \tau) \nonumber \\
    &&
  \end{eqnarray}

  \item{Equation ($N_0 -1$)C}
  
  \begin{flushleft}$y(r - (N_0 -2), 0 ,t_2, \ldots ,t_{N_0-1}, t_{N_0})$\end{flushleft}
  \begin{eqnarray}
    &=& \hat{y}(r - (N_0 - 1), 
            \tau - t_{N_0}, (\tau - t_{N_0}) + t_2, \ldots , (\tau - t_{N_0}) + t_{N_0 -1}, \tau) \nonumber \\
    &&
  \end{eqnarray}

  \item{Equation ($N_0$)B}
  
  \begin{flushleft}$\hat{y}(r - (N_0 -1), t_1 ,t_2, \ldots, t_{N_0-1}, \tau)$\end{flushleft}
  \begin{eqnarray}
    &=& \int_{x=0}^{\tau - t_{N_0-1}}{} \lambda_{r-(N_0-2)} e^{-\lambda_{r-(N_0-2)} x} \nonumber \\
    && \qquad * y(r-(N_0-2)),0, x + t_1, x + t_{N_0 - 1}) dx  \nonumber \\
    && + \lambda_{r-(N_0-2)} e^{-\lambda_{r-(N_0-2)} (\tau - t_{N_0 - 1}) } \nonumber \\
    &&  \qquad * \hat{y}(r-(N_0-2),0,
            (\tau - t_{N_0 - 1}) + t_1, \ldots , (\tau - t_{N_0 - 1}) + t_{N_0 -2}, \tau) \nonumber \\
    &&
  \end{eqnarray}

\end{enumerate}

It can be verified that the solutions to the above non-normalized equations are

\begin{eqnarray}
  \bar{y}(r+1, \vec{0}^{N_0} ) &=& \frac{C}{ \prod_{i=r-(N_0-2)}^{r}{\lambda_i} } \\
  \hat{y}(l, \vec{t}, \tau ) &=& \frac{C}{ \prod_{i=r-(N_0-2)}^{l}{\lambda_i} } \nonumber \\
                             && \qquad \qquad \text{ for } r-(N_0-2) \le l \le r\\
  y(l, \vec{t}, \tau ) &=& \frac{C}{ \prod_{i=r-(N_0-2)}^{l-1}{\lambda_i} } \nonumber \\
                             && \qquad \qquad \text{ for } r-(N_0-3) \le l \le r\\
  y(r-(N_0 -2), \vec{t}) = C \\
  \hat{y}(r-(N_0 -1), \vec{t}, \tau) = C
\end{eqnarray}

\emph{Remark:} These equilibrium solutions do not depend on the times remaining on outstanding orders.

\subsection{Expected time per $M$ state spent at inventory levels $l$}

From the appendix, we have $h(l,l,t) = \frac{1}{\lambda_{l}}(1- e^{-\lambda_{l} t} )$. Using this, and the derivations of $c(s,l)$ we find

\begin{equation}
    c((l',t_1,t_2, \ldots, t_{N_0}),l) =
      \left\{
      \begin{array}{rl}
        \frac{1}{\lambda_l}  & \text{ if } r <l \le l'\\
         \frac{1}{\lambda_{l}}(1- e^{-\lambda_{l} t_{N_0 - k +1} } )  & \text{ if } l' = l = r-(k-1) \text{ and } 1 \le k \le N_0-1 \\
         t_1  & \text{ if } r-(N_0-1) = l = l'\\
          0                    & \text{ otherwise }
      \end{array}
      \right.
\end{equation}

Substitution of the formulas for $c(s,l)$, $\bar{y}()$, $y()$ and $\hat{y}()$ into the equations for $w(l)$ found earlier, followed by direct integration and simplification yields

\begin{eqnarray}
  w(r+1) &=& \frac{C}{ \prod_{i=r-(N_0-2)}^{r+1}{\lambda_i} } \\
  w(l) &=& \frac{C \tau^{r-l+1}}{ (r-l+1)! \prod_{i=r-(N_0-2)}^{l}{\lambda_i} } \nonumber \\
        && \qquad \qquad \text{ for } r-(N_0-3) \le l \le r \\
  w(r - (N_0 -2)) &=& \frac{C \tau^{N_0 -1}}{ (N_0 -1)! \lambda_{r - (N_0 -2)} } \\
  w(r - (N_0 -1)) &=& \frac{C \tau^{N_0}}{ (N_0)! } \\
\end{eqnarray}

With $\beta_i = \lambda_i$ for $i>r-(N_0-1)$ and $\beta_{r-(N_0-1)} = 1$ these equations can be re-written as

\begin{equation}
   w(l) = C \frac{ \frac{1}{(r-l+1)!} \prod_{i=l+1}^{r+1}{ ( \tau \beta_i ) } } { \prod_{i = r-(N_0-1)}^{r+1}{\beta_i} } 
\end{equation}

Finally, we have

\begin{eqnarray}
 a(l) &=& \displaystyle{ \frac{w(l)}{ \sum_{i=l_L}^{r+q}{w(i)} } } \\
     &=& \displaystyle{ 
            \frac{ 
                    \displaystyle{ 
                       \frac{
                              \frac{1}{(r-l+1)!} \prod_{i=l+1}^{r+1}{ ( \tau \beta_i ) } } { \prod_{i = r-(N_0-1)}^{r+1}{\beta_i} } 
                     }
                  }
                  {
                    \sum_{j=r-(N_0-1)}^{r+1} 
                      \left\{  
                    \displaystyle{ 
                                    \frac{ \frac{1}{(r-j+1)!} \prod_{i=j+1}^{r+1}{ ( \tau \beta_i ) } } { \prod_{i = r-(N_0-1)}^{r+1}{\beta_i} } 
                                  }
                     \right\}
                  }
           } \\
      &=& \displaystyle{ 
            \frac{ 
                    \displaystyle{ 
                        \frac{1}{(r-l+1)!} \prod_{i=l+1}^{r+1}{ ( \tau \beta_i ) }  
                     }
                  }
                  {
                    \sum_{j=r-(N_0-1)}^{r+1} 
                      \left\{  
                    \displaystyle{ 
                                     \frac{1}{(r-j+1)!} \prod_{i=j+1}^{r+1}{ ( \tau \beta_i ) } 
                                  }
                     \right\}
                  }
           }
\end{eqnarray}

We notice that the equations do not depend on $\beta_{r-(N_0-1)}$. Therefore we can write the final formula in terms of $\lambda_i$ :

\begin{equation}
  a(l) = \frac{ 
                \displaystyle {\frac{1}{(r-l+1)!} }
                        \displaystyle
                        {
                            \prod_{i=l+1}^{r+1}{ ( \tau \lambda_i ) }  
                        }
               }
               {
                \displaystyle
                {
                  \sum_{j=r-(N_0-1)}^{r+1} 
                    {
                      \left\{  
                         \frac{1}{(r-j+1)!} \prod_{i=j+1}^{r+1}{ ( \tau \lambda_i ) } 
                       \right\}
                      }
                }
           } 
\end{equation}
\begin{flushright} $\qed$ \end{flushright}

\section{Application of the case $q=1$}

We consider a system of two stores with continuous review buy-as-sold policies $(r_1,1$ and $(r_2,1)$. Each store has a lateral transshipment cutoff level $c_n$. If a store's inventory level is at least $c_n$, outgoing lateral transshipments are allowed to the other store. The lateral transshipments take no time. The stores face customer demand for one unit at a time, which is Poisson distributed with rates $\gamma_n$. If a customer demand arrives when a store has no on-hand inventory, then a transshipment from the other store is made if allowed, otherwise the sale is lost.

Let $\beta_i$ represent the fraction of time store $i$ has no on-hand inventory ( $l=0$ ). The \emph{effective rates} at which inventory level decreases are given by ( for $m \ne n$ )

\begin{equation}
    \lambda_l^m =
      \left\{
      \begin{array}{rl}
        \gamma_m + \beta_n \gamma_n & \text{ if } l \ge c_m\\
        \gamma_m  & \text{ if } 1 \le l < c_m \\
        0 & \text{ if } l=0 
      \end{array}
      \right.
\end{equation}

We can compute $\beta_m = a(0)_m$ from the formulas for the average fraction of time spent at inventory level $0$.

\begin{eqnarray}
  \beta_m &=& a(0)_m \\
          &=& \frac{ 
                    \displaystyle {\frac{1}{(r_m+1)!} }
                            \displaystyle
                            {
                                \prod_{i=1}^{r_m+1}{ ( \tau \lambda_m ) }  
                            }
                   }
                   {
                    \displaystyle
                    {
                      \sum_{j=0}^{r+1} 
                        {
                          \left\{  
                             \frac{1}{(r_m-j+1)!} \prod_{i=j+1}^{r_m+1}{ ( \tau \lambda_m ) } 
                           \right\}
                          }
                    }
               } \\
    \beta_m &=& \frac{ 
                    \displaystyle {\frac{1}{(r_m+1)!} }
                            \displaystyle
                            {
                                \prod_{i=1}^{c_m -1}{ ( \tau \lambda_m ) }  
                                \prod_{i=c_m}^{r_m+1}{ \tau ( \lambda_m  + \beta_n \lambda_n) } 
                            }
                   }
                   {
                    \displaystyle
                    {
                      \sum_{j=0}^{r+1} 
                        {
                          \left\{  
                             \frac{1}{(r_m-j+1)!} 
                             \prod_{i=j+1}^{c_m -1}{ ( \tau \lambda_m ) } 
                             \prod_{i=max(c_m,j+1)}^{r_m+1}{  \tau ( \lambda_m  + \beta_n \lambda_n) } 
                           \right\}
                          }
                    }
               } \nonumber \\
            && 
\end{eqnarray}

with $m \ne n$ and the convention that an empty product equals $1$. This gives us two equations in two unknowns. We may solve for $\beta_1$ and $\beta_2$ (numerically or otherwise) and then compute $a(l)_i$, from which we can obtain the costs (holding costs, lost sales penalties, etc.) which are of interest for an inventory system.

\section{The case when $r=2$ and $q=2$}

We will solve the embedded equilibrium equations for the case when $r=2$ and $q=2$, where $\lambda_0 = 0$ and $\lambda_i=\lambda$ for $i \ge 1$. The embedded equilibrium equations are tractable in this case. We will be satisfied to solve them up to a constant multiple. To ease calculations, we will choose $\hat{y}(2,0,\tau) =1$. We also note that in this case we have (since $\lambda_0 = 0$ ) $g(0,0,\tau-t)=1$. Therefore $y(2,0,t) = \hat{y}(0, \tau-t,\tau)$ and $y(1,0,t) = 0$. In addition, from the appendix we have $f(2,1,x) = \lambda^2 x e^{-\lambda x}$. 

Rewrite the equation for $\hat{y}(0,t,\tau)$ as 

\begin{eqnarray}
  \hat{y}(0,t,\tau) 
      &=& \int_{x=0}^{\tau-t}{ f(2,1,x)y(2,0,t+x) dx} + f(2,1,\tau-t)\hat{y}(2,0,\tau) \nonumber \\
      &=& \int_{x=0}^{\tau-t}{ \lambda^2 x e^{-\lambda x} \hat{y}(0, \tau-t-x,\tau) dx} 
            + \lambda^2 (\tau-t) e^{-\lambda (\tau-t)}  \nonumber \\
  z(t) &=& \int_{x=0}^{\tau-t}{ \lambda^2 (\tau-t-x) e^{-\lambda (\tau-t-x)} z(x) dx} 
            + \lambda^2 (\tau-t) e^{-\lambda (\tau-t)} \nonumber \\
\end{eqnarray}

where $z(t) = \hat{y}(0,t,\tau)$. By repeatedly differentiating $z(t)$ and substituting $t=\tau-t$ this equation can be reduced to the ordinary differential equation

\begin{equation}
  0=-2\,{\lambda}^2 \frac{d^2}{d{t}^2} z(t) +\frac{d^4}{d{t}^4} z( t) 
\end{equation}

The general solution is $z(t) = A + Bt + C e^{\lambda t \sqrt{2}} + D e^{- \lambda t \sqrt{2}}$. Substituting this and solving the resulting equations we can obtain:

\begin{eqnarray}
A &=& \frac{1}{2}{\frac {\lambda \sqrt{2}( -2 \sqrt{2}+3+{ e^{\lambda \tau \sqrt{2}}} ) }{{ e^{\lambda\tau \sqrt{2}}}-3+2 \sqrt {2}}} \\
B &=& 0 \\
C &=& -\frac{1}{2} \frac{\lambda \sqrt{2}}{{ e^{ \lambda \tau \sqrt{2}}}-3+2 \sqrt {2}} \\
D &=& -\frac{1}{2} \frac{ ( -4+3 \sqrt{2} ) \lambda{ e^{\lambda \tau \sqrt{2}}}}{ e^{ \lambda \tau \sqrt{2}}-3+2 \sqrt{2}}
\end{eqnarray}

The values $a(l)$ can be obtained by direct integration following the procedure in section 6. The resulting formulas are complicated and we won't display them here.

\section{Appendix}

We review some useful results concerning sums of exponential random variables. Throughout this section $L$ refers to the Laplace transform operator with respect to time $t$.

\subsection{Probability density of a sum of exponential random variables}

Let $\{X_i\}$ be a collection of exponential random variables with rates ${\lambda_i}$. 

\begin{definition}
Let $X_{l_1,l_2} = \sum_{i=l_1}^{l_2}{X_i}$ with $l_1 \le l_2$.
\end{definition}

\begin{definition}
Define $f(l_2,l_1,t)$ to be the probablility density function of $X_{l_1,l_2}$. Define $\tilde{f}(l_2,l_1,s)$ to be $L\{f(l_2,l_1,t)\}(s)$, where $L$ is the Laplace transform operator with respect to $t$.
\end{definition}

The density function of a sum of variables is a convolution, so the Laplace transform becomes

\begin{eqnarray}
\tilde{f}(l_2,l_1,s) &=& \prod_{i=l_1}^{l_2}{\frac{\lambda_i}{\lambda_i + s}}
\end{eqnarray}

When the $\lambda_i$ are identical, $f(l_2,l_1,t)$ is an Erlang density. When all $\lambda_i$ are distinct, $f(l_2,l_1,t)$ is given by

\begin{equation}
f(l_2,l_1,t) =
    \left\{
    \begin{array}{rl}
      \lambda_{l_1}e^{-\lambda_{l_1}t} &\text{if } l_1=l_2 \\
      \displaystyle{ 
                      \sum_{i=l_1}^{l_2}{ \prod_{ \substack{j=l_1 \\ j \ne i} }^{l_2}{
                      \left( \frac{\lambda_j}{\lambda_j-\lambda_i} \right) \lambda_i e^{-\lambda_i t}  } }
                   } &\text{if } l_1 > l_2 \\
    \end{array}
    \right.
\end{equation}

A closed form for $f(l_2,l_1,t)$ is known for an arbitrary collection of $\lambda_i$ \cite{AM1}.

\subsection{Probability of decrease in inventory level before an order arrives}

\begin{definition}
Suppose the inventory level is currently $l'$ and an order is due to arrive in time $t$ from now, with no orders arriving earlier. Define $g(l,d,t)$ to be the probability that the inventory level decreases from $l$ to $l-d$ in time $t$. Define $\tilde{g}(l,d,s) = L\{g(l,d,t)\}(s)$.
\end{definition}

We have

\begin{equation}
  g(l,d,t) = 
      \left\{
      \begin{array}{rl}
        \displaystyle{ \int_{x=t}^{\infty} {\lambda_{l}  e^{-\lambda_{l} x} dx } } & \text{if } d = 0 \\
        \displaystyle{ \int_{x=0}^{t} { f(l,l-d+1,x) e^{-\lambda_{l-d}(t-x)} dx }} & \text{if } d \ge 1
      \end{array}
      \right.
\end{equation}

and so using the convolution property of the Laplace transform

\begin{equation}
  \tilde{g}(l,d,s) = 
    \left\{
    \begin{array}{rl}
      \displaystyle{ \frac{1}{\lambda_{l} + s} }
            & \text{if } d = 0 \\
      \displaystyle{ \frac{1}{\lambda_{l-d} + s} \prod_{i=l-d+1}^{l}{\frac{\lambda_i}{\lambda_i +s}} } 
            & \text{if } d \ge 1
    \end{array}
    \right.
\end{equation}

If all $\lambda_i$ are distinct, direct integration of $ \int_{x=0}^{t} { f(l,l-d+1,x) e^{-\lambda_{l-d}(t-x)} dx } $ gives the result
\\
\\
$  g(l,d,t) = $
\\
\begin{equation}
    \left\{
    \begin{array}{rl}
      \displaystyle{ e^{-\lambda_l t} }
            & \text{if } d = 0 \\
      \displaystyle{ \left( \frac{\lambda_l}{\lambda_{l-1} - \lambda_l}\right) \left( e^{-\lambda_l t} - e^{-\lambda_{l-1} t} \right) } 
            & \text{if } d = 1 \\
      \displaystyle{ 
                      \sum_{i=l-d+1}^{l}
                        { \left[ \prod_{\substack{j=l-d+1 \\ j \neq i}}^{l}
                           {
                             \left( \frac{\lambda_j}{\lambda_j - \lambda_i} \right) 
                             \left(  \frac{\lambda_i}{\lambda_{l-d} - \lambda_i} \right)      
                             \left( e^{-\lambda_i t} - e^{-\lambda_{l-d} t} \right)
                           } 
                        \right] }
                    }
            & \text{if } d > 1
    \end{array}
    \right.
\end{equation}

\subsection{Expected time spent in an inventory level before an order arrives}

\begin{definition}
Suppose the inventory level is currently $l'$ and an order is due to arrive in time $t$ from now, with no orders arriving earlier. Let $h(l',l,t)$ be the expected time spent in inventory level $l, l \le l'$ during time $t$. Let $\tilde{h}(l',l,s) = L\{h(l',l,t)\}(s) $.
\end{definition}

A time $x$ is spent in level $l$ when either A) $X_{l',l+1} \le t-x$ and $X_{l,l} = x$ or B) $X_{l',l+1} = t-x$ and $X_{l,l} \ge x$. Integrating these contributons gives us

\begin{eqnarray}
h(l',l,t) = & \displaystyle{ \int_{0}^{t}{ \left(x \lambda_l e^{-\lambda_l x} \right)
                             \left( \int_{0}^{t-x}{f(l',l+1,s) ds} \right) dx }
                           } \nonumber \\
           & + \displaystyle{ \int_{0}^{t}{ f(l',l+1,t-x)
                             \left( x \int_{x}^{\infty}{ \lambda_l e^{-\lambda_l s} ds} \right) dx }
                           } \\
          = & h_1(l',l,t) + h_2(l',l,t)
\end{eqnarray}

$h_1$ and $h_2$ are convolutions, so we take the Laplace transforms of $h_1$ and $h_2$ to get

\begin{eqnarray}
\tilde{h_1}(s) =& \displaystyle{ \frac{1}{s} \left( \frac{\lambda_l}{ (\lambda_l +s)^2 } \right)
                                  \prod_{i=l+1}^{l'}{ \frac{\lambda_i}{\lambda_i +s} }
                                } \\
\tilde{h_2}(s) =& \displaystyle{  \frac{1}{ (\lambda_l +s)^2 }
                                  \prod_{i=l+1}^{l'}{ \frac{\lambda_i}{\lambda_i +s} }
                                } 
\end{eqnarray}

Summing gives the Laplace transform of $h$

\begin{eqnarray}
\tilde{h}(l',l,s) =& \displaystyle{ \frac{1}{ (\lambda_l +s)^2 } \left( \frac{\lambda_l}{s} +1  \right) 
                                  \prod_{i=l+1}^{l'}{ \frac{\lambda_i}{\lambda_i +s} }
                                } \\
                  =& \displaystyle{ \frac{1}{ s(\lambda_l +s) }
                                  \prod_{i=l+1}^{l'}{ \frac{\lambda_i}{\lambda_i +s} }
                                } 
\end{eqnarray}

It can be verified either by applying the inverse transform $L^{-1}$ to $\tilde{h}$ or by integrating the equation for $h(l',l,t)$ directly that (for distinct $\lambda_i$) $h$ has the form
\\
\\
$h(l',l,t) = $
\begin{equation}
  t  \text{ , if } l = l' \text{ and } \lambda_{l'} = 0
\end{equation}
\begin{equation}
  \frac{1}{\lambda_{l'}}(1- e^{-\lambda_{l'} t} )
     \text{ , if } l = l' \text{ and } \lambda_{l'} > 0
\end{equation}
\begin{equation}
  \frac{1}{\lambda_{l'-1}}(1- e^{-\lambda_{l'-1} t} )
      + ( \frac{1}{\lambda_{l'-1}-\lambda_{l'}} ) ( e^{-\lambda_{l'-1} t} - e^{-\lambda_{l'} t} 
     \text{ , if } l = l' - 1
\end{equation}
\begin{equation}
  \sum_{i=l+1}^{l'} 
    { \prod_{\substack{j=l+1 \\ j \neq i}}^{l'} 
      { 
       ( \frac{\lambda_j}{\lambda_j - \lambda_i} )
        \left\{
          \frac{1}{\lambda_l} (1-e^{-\lambda_l t} )
            + ( \frac{1}{\lambda_i-\lambda_l} ) ( e^{-\lambda_i t} - e^{-\lambda_l t}  )
        \right\}
      } 
    }
   \text{ , if } l < l' - 1
\end{equation}

\begin{remark}
Direct numerical calculation using the formulas $f$, $g$ and $h$ tends to explode for distinct $\lambda_i$ when some $\lambda_i$ are close together. Numerical inversion of the Laplace transforms $\tilde{f}$, $\tilde{g}$ and $\tilde{h}$ avoids this problem.
\end{remark}

\end{document}